\documentclass{article}
\usepackage{amsfonts}
\usepackage{amsmath}
\usepackage{amsmath}
\usepackage{amssymb}

\newtheorem{theorem}{Theorem}

\newtheorem{corollary}[theorem]{Corollary}

\newtheorem{lemma}[theorem]{Lemma}

\begin{document}

\title{Extinction of decomposable branching processes\thanks{This work is supported by the RSF under a grant 14-50-00005. } }
\author{Vatutin V.A.\thanks{%
Department of Discrete Mathematics, Steklov Mathematical Institute, 8,
Gubkin str., 119991, Moscow, Russia; e-mail: vatutin@mi.ras.ru}, Dyakonova
E.E.\thanks{%
Department of Discrete Mathematics, Steklov Mathematical Institute, 8,
Gubkin str., 119991, Moscow, Russia; e-mail: elena@mi.ras.ru}}
\date{}
\maketitle

\date{}
\maketitle

\begin{abstract}
The asymptotic behavior, as $n\rightarrow \infty $ of the conditional distribution of the number of particles in a decomposable
critical branching process $\mathbf{Z}%
(m)=(Z_{1}(m),...,Z_{N}(m)),$  with $N$ types of particles
 at moment $m=n-k,\, k=o(n),$ is investigated given that the extinction moment of the process is $n$.
\end{abstract}



\section{Introduction}
We consider a Galton-Watson branching process with $N$ types of particles
labelled $1,2,...,N$ in which a type
~$i$ parent particle may produce children of types $j\geq i$ only. Let $\eta _{i,j}$ be the number of type $j$ children produced by a type~$i$ parent particle. According to our assumption $\eta _{i,j}=0$ if $i>j.$

 In what follows we rather often use the offspring generating function of type $N$ particles. For this reason, to simplify notation we put
\begin{equation}
h(s)=\mathbf{E}\left[
\,s^{\eta _{N,N}}\right].  \label{Simple}
\end{equation}

Denote by $\mathbf{e}_{i}$ the $N$-dimensional vector whose $i$-th component is equal to one while the remaining are equal to zero and let
 $\mathbf{0}=(0,...,0)$ be an $N$-dimensional vector all whose components are equal to zero.

Let
\begin{equation*}
\mathbf{Z}(n)=(Z_{1}(n),...,Z_{N}(n))
\end{equation*}
be the population size at moment $n\in \mathbb{Z}_{+}=\left\{ 0,1,...\right\}.$
We denote by \begin{equation*}
m_{i,j}(n)=\mathbf{E}\left[ Z_{j}(n)|\mathbf{Z}\left( 0\right)
=\mathbf{e}_{i}\right]
\end{equation*}
the expectations of the components of the vector $\mathbf{Z}(n).$  Let  $$m_{i,j}=m_{i,j}(1)=\mathbf{E}[\eta_{i,j}]$$
be the average number of the direct descendants of type~$j$ generated by a type~$i$ particle.

We say that \textbf{Hypothesis A} is valid if the decomposable branching
process with $N$ types of particles is strongly critical, i.e. (see \cite%
{FN2})
\begin{equation}
m_{i,i}=\mathbf{E}\left[ \eta _{i,i}\right] =1,\quad i=1,2,...,N
\label{Matpos}
\end{equation}%
and, in addition,
\begin{equation}
m_{i,i+1}=\mathbf{E}\left[ \eta _{i,i+1}\right] \in (0,\infty ),\
i=1,2,...,N-1,  \label{Maseq}
\end{equation}%
\begin{equation}
\mathbf{E}\left[ \eta _{i,j}\eta _{i,k}\right] <\infty ,\,i=1,...,N;\
k,j=i,i+1,...,N,  \label{FinCovar}
\end{equation}%
with
\begin{equation}
b_{i}=\frac{1}{2}Var\left[ \eta _{i,i}\right] \in \left( 0,\infty \right)
,i=1,2,...,N.  \label{FinVar}
\end{equation}%
Thus, a particle of the process is able to produce the direct descendants of
its own type, of the next in the order type, and (not necessarily, as direct
descendants) of all the remaining in the order types, but not any preceding
ones.

In the sequel we assume (if otherwise is not stated) that $\mathbf{Z}(0)=%
\mathbf{e}_{1}$, i.e. we suppose that the branching process in question is
initiated at time $n=0$ by a single particle of type 1.

Denote by $T_{N}$ the extinction
moment of the process. The aim of the present paper is to investigate the
asymptotic behavior, as $n\rightarrow \infty $ of the conditional distribution of the number of particles
in the decomposable critical  branching process $\mathbf{Z}(m)=(Z_{1}(m),...,Z_{N}(m))$ with $N$  types of particles
at moment $m=n-k,k=o(n)$  given $T_{N}=n.$

 Decomposable branching processes of different structure have been investigated in a number of papers.

We mention in this connection paper \cite{Afa15} dealing with the structure of the two-type decomposable critical Galton-Watson branching processes
 in which the total number of type 1 particles is fixed,  and articles  \cite{FN}--\cite{Z} in which, for  $N$-type decomposable critical Markov branching processes asymptotic representations for the probability of the event  $\left\{ T_{N}>n\right\} $ are found and Yaglom-type limit theorems are proved describing (under various restrictions) the distribution of the number of particles in these processes (and their reduced analogues) under the condition $T_{N}>n$.

In \cite{VD15} the decomposable critical  branching processes obeying the conditions of Hypothesis $A$ were investigated under the assumption $T_{N}=n.$  In the mentioned paper the conditional limit distributions are found  for the properly scaled components of the vector
 \begin{equation*}
\mathbf{Z}(m)=(Z_{1}(m),...,Z_{N}(m))
\end{equation*}
given that the parameter $m=m(n)$ varies in such a way that
\begin{equation*}
\lim_{n\rightarrow \infty }\frac{m}{n}=x\in \lbrack 0,1).
\end{equation*}

The following statement established in \cite{VD15} is of a particular interest:
\begin{theorem}
\label{T_finalstage}If Hypothesis $A$ is valid and $m=m(n)\to\infty$ as $n\to\infty$ in such a way that $m\sim xn,x\in \left(
0,1\right) ,$ then, for any $s_{i}\in \lbrack 0,1],i=1,2,...,N-1$ and $\lambda _{N}>0$
\begin{eqnarray*}
&&\lim_{n\rightarrow \infty }\mathbf{E}\left[ s_{1}^{Z_{1}(m)}\cdot \cdot
\cdot s_{N-1}^{Z_{N-1}(m)}\exp \left\{ -\lambda
_{N}\frac{Z_{N}(m)}{b_{N}n}\right\} \Big|\,T_{N}=n\right] \\
&&\qquad \qquad \qquad =\left( \frac{1+\lambda _{N}(1-x)}{1+\lambda
_{N}x\left( 1-x\right) }\right) ^{-1+1/2^{N-1}}\frac{1}{\left( 1+\lambda
_{N}x\left( 1-x\right) \right) ^{2}}.
\end{eqnarray*}
\end{theorem}

It follows from Theorem \ref{T_finalstage} that if $T_{N}=n$ and the parameter $m=m(n)$ varies within the specified range,
then the population consists (in the limit) of type $N$ particles only.

In the present paper,   complimenting paper \cite{VD15}, we concentrate on the case
\begin{equation*}
\lim_{n\rightarrow \infty }\frac{m}{n}=1.
\end{equation*}

\begin{theorem}
\label{T_death} If Hypothesis $A$ is valid and $k=k(n)=$ $n-m\rightarrow
\infty $ as $n\rightarrow \infty $ in such a way that $k=o(n),$ then
 \begin{equation*}
\lim_{n\rightarrow \infty }\mathbf{E}\left[ s_{1}^{Z_{1}(m)}\cdot \cdot
\cdot s_{N-1}^{Z_{N-1}(m)}\exp \left\{ -\lambda
_{N}\frac{Z_{N}(m)}{b_{N}k}\right\} \Big|\,T_{N}=n\right] =\frac{1}{\left( 1+\lambda _{N}\right) ^{2}}.
\end{equation*}
\end{theorem}

Denote by $h_{n}(s)$ the $n-$th iteration of the probability generating function $h(s)$.
It is known (see, for instance, \cite{AN},
page 93, formula (16)), that
 \begin{equation}
\lim_{n\rightarrow \infty }b_{N}n^2\left( h_{n}(s)-h_{n}(0)\right) =U(s)
\label{DefHarm}
\end{equation}
exists for any fixed $s\in \lbrack 0,1),$ where где $U(s)$ is the generating function of the so-called harmonic measure.
In addition,
\begin{equation}
U(h(s))=U(s)+1, \quad s\in \lbrack 0,1). \label{harm}
\end{equation}

The following theorem, complimenting Theorem \ref{T_death} describes the final stage of the development of the process
 given its extinction at a distant moment $n$.

\begin{theorem}
\label{T_deathFin} If Hypothesis $A$ is valid and $k=n-m=const$
as $n\rightarrow \infty ,$ then
 \begin{equation*}
\lim_{n\rightarrow \infty }\mathbf{E}\left[ s_{1}^{Z_{1}(m)}\cdot \cdot
\cdot s_{N-1}^{Z_{N-1}(m)}s_{N}^{Z_{N}(m)}\Big|\,T_{N}=n\right] =s_{N}\left(
U\left( s_{N}h_{k+1}(0)\right) -U\left( s_{N}h_{k}(0)\right)
\right).
\end{equation*}
\end{theorem}

\textbf{Remark}. We know by (\ref{harm}) that
$$
U\left(h_{k+1}(0)\right) -U\left(h_{k}(0)\right)=1
$$
and, therefore, the limit distribution we have found in Theorem \ref{T_deathFin} is proper.

\section{\protect\bigskip Auxiliary results}

We use the symbols $\mathbf{P}_{i}$ and $\mathbf{E}_{i}$ to denote the
probability and expectation calculated under the condition that a branching
process is initiated at moment $n=0$ by a single particle of type $i$.
Sometimes we write $\mathbf{P}$ and $\mathbf{E}$ for $\mathbf{P}_{1}$ and $%
\mathbf{E}_{1},$ respectively.

Introduce the constants
\begin{equation}
 c_{N,N}=1/b_{N},\quad c_{i,N}=\left( \frac{1}{b_{N}}\right) ^{1/2^{N-i}}\prod_{j=i}^{N-1}\left(
\frac{m_{j,j+1}}{b_{j}}\right) ^{1/2^{j-i+1}},\quad i<N.
\label{Const2}
\end{equation}
and
\begin{equation}
D_{i}=(b_{i}m_{i,i+1})^{1/2^{i}}c_{1,i},i=1,2,...,N.  \label{DcConnection}
\end{equation}
It is not difficult to check that
\begin{equation}
c_{1,N}=D_{N-1}\left( \frac{1}{b_{N}}\right) ^{1/2^{N-1}}=D_{N-1}\left(c_{N,N}\right) ^{1/2^{N-1}}.  \label{Const3}
\end{equation}
Denote
\begin{equation*}
T_{ki}=\min \left\{ n\geq
1:Z_{k}(n)+Z_{k+1}(n)+...+Z_{i}(n)=0|\mathbf{Z}(0)=\mathbf{e}_{k}\right\}
\end{equation*}
the extinction moment of the population consisting of particles of types
 $k,k+1,...,i,$ given that the process was initiated at time $n=0$ by a single particle of type $k.$ To simplify formulas we set $T_{i}=T_{1i}.$

We fix
$N\geq 2$ and use, when it is needed, the notation
\begin{equation*}
\gamma _{0}=0,\ \gamma _{i}=\gamma _{i}(N)=2^{-(N-i)},\ i=1,2,...,N.
\end{equation*}

The starting point of our arguments is the following theorem proved in  \cite{FN} (see also \cite{FN2}):

\begin{theorem}
\label{T_Foster}Let $\mathbf{Z}(n),n=0,1,...,$
be a decomposable branching process meeting  conditions
(\ref{Matpos}), (\ref{Maseq}) and (\ref{FinCovar}).
Then, as $n\rightarrow \infty $
\begin{equation}
\mathbf{P}_{i}(\mathbf{Z}(n)\neq
\mathbf{0})\sim c_{i,N}n^{-1/2^{N-i}},  \label{SurvivSingle}
\end{equation}
where $c_{i,N}$ are the same as in (\ref{Const2}).
\end{theorem}

This result was complemented in \cite{VD15} by the following two statements the first of which is a local limit theorem.

\begin{theorem}
\label{T_loc}(see \cite{VD15}) If Hypothesis $A$ is valid, then, as  $n\rightarrow \infty $
\begin{equation}
\mathbf{P}\left( T_{iN}=n\right) \sim \frac{g_{i,N}}{n^{1+\gamma
_{i}}},\quad i=1,2,...,N,  \label{LocAsymp}
\end{equation}
where
\begin{equation}
g_{i,N}=\gamma _{i}c_{i,N}.  \label{Const4}
\end{equation}
\end{theorem}

\begin{corollary}
(see \cite{VD15}) If $\sqrt{n}\ll l\ll
n$, then
\begin{equation}
\lim_{n\rightarrow \infty }\mathbf{P}(Z_{1}(l)+\cdots
+Z_{N-1}(l)>0|T_{N}=n)=0.  \label{NoPrevious}
\end{equation}
\end{corollary}

Let $\eta _{r,j}\left( k,l\right) $ be the number of type $j$ daughter
particles of the $l-$th particle of type~$r$ belonging to the $k-$th
generation and let
 \begin{equation*}
W_N=\sum_{r=1}^{N-1}\sum_{k=0}^{T_{r}}\sum_{q=1}^{Z_{r}(k)}\eta _{r,N}\left(
k,q\right)
\end{equation*}
be the total number of type $N$ daughter particles generated by all
the particles of types $1,2,...,N-1$ ever born in the process given that the
process is initiated at time $n=0$ by a single particle of type $1.$

Asymptotic properties of the tail distribution of the random variable $W_N$ are described in the next lemma.

\begin{lemma}
\label{L_Laplace} (see \cite{VV14}, Lemma 1). Let Hypothesis $A$ be valid. Then, as $\theta \downarrow 0$
\begin{equation}
1-\mathbf{E}\left[ e^{-\theta W_{N}}\,|\mathbf{Z}(0)=\mathbf{e}_{1}\right]
\sim D_{N-1}\theta ^{\gamma _{1}}.  \label{Tot1}
\end{equation}
\end{lemma}

Proving the main results of the paper we will relay on the following statement.
\begin{lemma}
\label{L_Diff1}If  $m_{N,N}=1,\, b_{N}\in \left( 0,\infty
\right) $ and $k=k(n)=$ $n-m\rightarrow \infty $ as
$n\rightarrow \infty $ in such a way that $k=o(n),$ and
\begin{equation*}
s=\exp \left\{ -\frac{\lambda }{b_{N}k}\right\}, \lambda>0,
\end{equation*}
then
 \begin{equation*}
\lim_{n\rightarrow \infty }\frac{b_{N}\lambda n^{2}}{k}\left(
h_{m}(s)-h_{m}(0)\right) =1.
\end{equation*}
\end{lemma}

\textbf{Proof}. Let $q=q(k)$ be a positive integer such that
\begin{equation*}
h_{q}(0)\leq s\leq h_{q+1}(0).
\end{equation*}
Then
\begin{equation*}
h_{m+q}(0)-h_{m}(0)\leq h_{m}(s)-h_{m}(0)\leq h_{m+q+1}(0)-h_{m}(0).
\end{equation*}
Clearly, $q(k)\rightarrow \infty $
as $k\rightarrow \infty $, and in view of the representation
\begin{equation*}
1-h_{q}(0)\sim \left( b_{N}q\right) ^{-1},\quad 1-s\sim \lambda \left(
b_{N}k\right) ^{-1}
\end{equation*}
and \begin{equation*}
\lim_{h\rightarrow \infty }\frac{1-h_{q}(0)}{1-h_{q+1}(0)}=1,
\end{equation*}
we have $
q\sim k\lambda^{-1},\, k\rightarrow \infty .$
This and the local limit theorem for the critical Galton-Watson processes  (see, for instance, \cite{AN}, Corollary 1.9.I, p. 23)
imply that, as $1\ll k\ll m\sim
n\rightarrow \infty $
\begin{eqnarray*}
h_{m}(s)-h_{m}(0) &\leq &h_{m+q+1}(0)-h_{m}(0)=\sum_{j=1}^{q}\left( h_{m+j+1}(0)-h_{m+j}(0)\right)  \\
&\sim &\sum_{j=1}^{q}\frac{1}{b_{N}\left( m+j\right) ^{2}}\sim
\frac{q}{b_{N}n^{2}}\sim \frac{k}{b_{N}\lambda n^{2}}.
\end{eqnarray*}
Thus,
\begin{equation*}
\limsup_{n\rightarrow \infty }\frac{b_{N}\lambda n^{2}}{k}\left(
h_{m}(s)-h_{m}(0)\right) \leq 1.
\end{equation*}
Similar arguments show that
\begin{equation*}
\liminf_{n\rightarrow \infty }\frac{b_{N}\lambda n^{2}}{k}\left(
h_{m}(s)-h_{m}(0)\right) \geq 1.
\end{equation*}

The lemma is proved.

\begin{lemma}
\label{L_Diff2}If Hypothesis $A$ is valid, then, for any $\lambda >0$
\begin{equation*}
\lim_{n\rightarrow \infty }\frac{1}{n^{1-\gamma _{1}}}\mathbf{E}\left[
W_{N}\exp \left\{ -\lambda\frac{W_{N}}{b_{N}n}\right\} \right] =\frac{\gamma
_{1}b_{N}c_{1,N}}{\lambda ^{1-\gamma _{1}}}=\frac{b_{N}g_{1,N}}{\lambda
^{1-\gamma _{1}}}.
\end{equation*}
\end{lemma}

\textbf{Proof}. Let
\begin{equation*}
W_{N}(k)=W_{N}^{(1)}+W_{N}^{(2)}+...+W_{N}^{(k)},
\end{equation*}
where the summands are independent random variables and $W_{N}^{(i)}\overset{d}{=}W_{N},i=1,...,k$.
Setting $\theta =\lambda \left( b_{N}n\right)
^{-1}$ in Lemma \ref{L_Laplace} and recalling (\ref{Const3}), we obtain
\begin{equation*}
\lim_{n\rightarrow \infty }n^{\gamma _{1}}\mathbf{E}\left[ 1-\exp \left\{
-\frac{\lambda }{b_{N}n}W_{N}\right\} \right] =D_{N-1}\left( \frac{\lambda
}{b_{N}}\right) ^{\gamma _{1}}=c_{1,N}\lambda ^{\gamma _{1}}.
\end{equation*}
Hence it follows that, for all $\lambda>0$
\begin{eqnarray}
&&\lim_{n\rightarrow \infty }\mathbf{E}\left[ \exp \left\{ -\frac{\lambda
}{b_{N}n}W_{N}(\left[ n^{\gamma _{1}}\right] )\right\} |\mathbf{Z}(0)=\left[
n^{\gamma _{1}}\right] \mathbf{e}_{1}\right]  \notag \\
&=&\lim_{n\rightarrow \infty }\mathbf{E}^{\left[ n^{\gamma _{1}}\right]
}\left[ \exp \left\{ -\frac{\lambda }{b_{N}n}W_{N}\right\}
|\mathbf{Z}(0)=\mathbf{e}_{1}\right] =\exp \left\{ -c_{1,N}\lambda ^{\gamma _{1}}\right\} .
\label{Anali1}
\end{eqnarray}
Since the sequence of functions under the limit consists of analytical and uniformly bounded functions in the domain  $\left\{
Re\lambda >0\right\} :$
\begin{equation*}
\left\vert \mathbf{E}^{\left[ n^{\gamma _{1}}\right] }\left[ \exp \left\{
-\frac{\lambda }{b_{N}n}W_{N}\right\} |\mathbf{Z}(0)=\mathbf{e}_{1}\right]
\right\vert \leq 1,
\end{equation*}
it follows from the Montel theorem (see \cite{Mar78}, Ch. VI, Section 7) that this sequence is compact. Moreover, since this sequence converges for real
$\lambda >0,$ the Vitali theorem (see \cite{Mar78}, Ch. VI, Section 8) and the uniqueness theorem for analytical functions imply convergence in (\ref{Anali1}) for all $\lambda$ satisfying the condition  $\left\{
Re\lambda >0\right\}.$
Moreover, according to the Weierstrass theorem  (see \cite{Mar78}, Ch. VI, Section 6)  the derivatives of the prelimiting functions converge to the derivative of the limiting function in the specified domain. Whence, on account of the equality
\begin{equation*}
\lim_{n\rightarrow \infty }\mathbf{E}\left[ \exp \left\{ -\frac{\lambda
}{b_{N}n}W_{N}\right\} |\mathbf{Z}(0)=\mathbf{e}_{1}\right] =1
\end{equation*}
it follows that
\begin{eqnarray*}
&&\lim_{n\rightarrow \infty }\frac{1}{b_{N}n^{1-\gamma _{1}}}\mathbf{E}\left[
W_{N}\exp \left\{ -\frac{\lambda }{b_{N}n}W_{N}
\right\} \right] \mathbf{E}^{\left[ n^{\gamma _{1}}\right] }\left[ \exp
\left\{ -\frac{\lambda }{b_{N}n}W_{N}\right\}
|\mathbf{Z}(0)=\mathbf{e}_{1}\right] \\
&&\qquad \qquad \qquad \qquad \qquad =-\frac{\partial }{\partial \lambda
}\exp \left\{ -c_{1,N}\lambda ^{\gamma _{1}}\right\}
\end{eqnarray*}
or, in view of (\ref{Anali1}) and (\ref{Const4})
 \begin{equation*}
\lim_{n\rightarrow \infty }\frac{1}{n^{1-\gamma _{1}}}\mathbf{E}\left[
W_{N}\exp \left\{ -\frac{\lambda }{b_{N}n}W_{N}\right\} \right] =\frac{\gamma _{1}b_{N}c_{1,N}}{\lambda ^{1-\gamma
_{1}}}=\frac{b_{N}g_{1,N}}{\lambda ^{1-\gamma _{1}}}.
\end{equation*}
The lemma is proved.

Let\begin{equation*}
I_{k}(m)=I\{Z_{1}(m)+\cdots +Z_{k}(m)=0\}
\end{equation*}
be the indicator of the event that there are no particles in of types $1,2,...,k$ in the population at time $m.$  We also agree to consider that $I_{0}(m)=1$.

\begin{corollary}
\label{C_Cansor}If Hypothesis $A$ is valid, then, for any $\lambda >0$
\begin{equation*}
\lim_{n\rightarrow \infty }\frac{1}{n^{1-\gamma _{1}}}\mathbf{E}\left[
W_{N}\exp \left\{ -\frac{\lambda }{b_{N}n}W_{N}\right\}
I_{N-1}(n^{2/3})\right] =\frac{b_{N}g_{1,N}}{\lambda ^{1-\gamma _{1}}}.
\end{equation*}
\end{corollary}

\textbf{Proof}. To check the validity of the statement of the lemma it is sufficient to note that, by virtue of (\ref{SurvivSingle})
 \begin{equation*}
\mathbf{P}\left( Z_{1}(n^{2/3})+\cdots +Z_{N-1}(n^{2/3})>0\right) =O\left(
\left( n^{2/3}\right) ^{-1/2^{N-2}}\right) =o\left( n^{-1/2^{N-1}}\right)
=o\left( n^{-\gamma _{1}}\right) ,
\end{equation*}
to make use of the equalities
\begin{eqnarray*}
&&\mathbf{E}\left[ \left( 1-\exp \left\{ -\frac{\lambda
}{b_{N}n}W_{N}\right\} \right) I_{N-1}(n^{2/3})\right] \\
&=&\mathbf{E}\left[ 1-\exp \left\{ -\frac{\lambda
}{b_{N}n}W_{N}I_{N-1}(n^{2/3})\right\} \right] \\
&=&\mathbf{E}\left[ 1-\exp \left\{ -\frac{\lambda }{b_{N}n}W_{N}\right\}
\right] -\mathbf{P}\left( Z_{1}(n^{2/3})+\cdots +Z_{N-1}(n^{2/3})>0\right) ,
\end{eqnarray*}
and, replacing $W_{N}$ by $W_{N}I_{N-1}(n^{2/3})$,
to repeat the arguments we have applied to prove Lemma \ref{L_Diff2}.

\begin{lemma}
\label{L_Diff3}If Hypothesis $A$ is valid, then, for any $\lambda >0$ and $1\ll k\ll n $
\begin{equation*}
\lim_{n\to\infty}\frac{n^{1+\gamma _{1}}}{k^{2}}\mathbf{E}\left[
Z_{N}(m)\exp \left\{ -\lambda \frac{Z_{N}(m)}{b_{N}k}\right\}
I_{N-1}(n^{2/3})\right] =\frac{b_Ng_{1,N}}{\lambda ^{2}}.
\end{equation*}
\end{lemma}

\textbf{Proof. } We consider the difference
\begin{equation*}
\Delta (m,k;\lambda )=\mathbf{E}\left[ \exp \left\{ -\lambda
\frac{Z_{N}(m)}{b_{N}k}\right\} I_{N-1}(n^{2/3})\right] -\mathbf{E}\left[ I\left\{
Z_{N}(m)=0\right\} I_{N-1}(n^{2/3})\right] .
\end{equation*}
Clearly,
\begin{equation*}
\frac{\partial \Delta (m,k;\lambda )}{\partial \lambda
}=-\frac{1}{b_{N}k}\mathbf{E}\left[ Z_{N}(m)\exp \left\{ -\lambda
\frac{Z_{N}(m)}{b_{N}k}\right\} I_{N-1}(n^{2/3})\right] .
\end{equation*}
Introduce, for $m\geq T_{N-1}$ the quantity
\begin{equation*}
\mathcal{H}_{m}\left( s\right)
=\prod_{r=1}^{N-1}\prod_{k=0}^{T_{r}}\prod_{l=1}^{Z_{r}(k)}\left(
h_{m-k}(s)\right) ^{\eta _{r,N}\left( k,l\right) }.
\end{equation*}
Then, for $n^{2/3}<m<n$
\begin{equation*}
\Delta (m,k;\lambda )=\mathbf{E}\left[ \left( \mathcal{H}_{m}\left(
s\right) -\mathcal{H}_{m}\left( 0\right) \right) I_{N-1}(n^{2/3})\right],
\end{equation*}
where
 \begin{equation}
s=\exp \left\{ -\frac{\lambda }{b_{N}k}\right\}. \label{defS}
\end{equation}

Observe that, by the criticality condition
\begin{equation*}
h_{l+1}(s)=h(h_{l}(s))\geq h_{l}(s)
\end{equation*}
and
\begin{equation*}
h_{l+1}(s)-h_{l+1}(0)=h(h_{l}(s))-h(h_{l}(0))\leq h_{l}(s)-h_{l}(0)
\end{equation*}
for any fixed
$s\in \left[ 0,1\right].$
This yelds
\begin{equation*}
\left( h_{m-T_{N-1}}(s)\right) ^{W_{N}}\leq \mathcal{H}_{m}\left( s\right)
\leq \left( h_{m}(s)\right) ^{W_{N}}.
\end{equation*}
It is not difficult to check that for any $s\in \left[
0,1\right] $ and $n^{2/3}<m<n$
\begin{eqnarray*}
&&\mathbf{E}\left[ s^{Z_{N}(m)}I_{N-1}(n^{2/3})\right] -\mathbf{E}\left[
I\left\{ Z_{N}(m)=0\right\} I_{N-1}(n^{2/3})\right]  \\
&=&\mathbf{E}\left[ \left( \mathcal{H}_{m}\left( s\right)
-\mathcal{H}_{m}\left( 0\right) \right) I_{N-1}(n^{2/3})\right] .
\end{eqnarray*}
Using the inequalities
\begin{eqnarray*}
\sum_{l=1}^{J}r_{l}\frac{\left( a_{l}-b_{l}\right)
}{b_{l}}b^{\sum_{j=1}^{J}r_{j}} &\leq &\sum_{l=1}^{J}r_{l}\frac{\left(
a_{l}-b_{l}\right) }{b_{l}}\prod\limits_{j=1}^{J}b_{j}^{r_{j}} \\
&\leq
&\prod\limits_{j=1}^{J}a_{j}^{r_{j}}-\prod\limits_{j=1}^{J}b_{j}^{r_{j}} \\
&\leq &\sum_{l=1}^{J}r_{l}\frac{\left( a_{l}-b_{l}\right)
}{a_{l}}\prod\limits_{j=1}^{J}a_{j}^{r_{j}}\leq \sum_{l=1}^{J}r_{l}\left(
a_{l}-b_{l}\right)a^{\sum_{j=1}^{J}r_{j}-1},
\end{eqnarray*}
being valid for $0< b\leq b_{j}\leq a_{j}\leq a\leq 1$
and nonnegative integers $r_{j},$ we obtain
\begin{eqnarray*}
&&W_{N}\left( h_{m-T_{N-1}}(0)\right) ^{W_{N}}\left(
h_{m}(s)-h_{m}(0)\right)  \\
&\leq &\sum_{r=1}^{N-1}\sum_{k=0}^{T_{r}}\sum_{l=1}^{Z_{r}(k)}\eta
_{r,N}\left( k,l\right) \frac{\left( h_{m-k}(s)-h_{m-k}(0)\right)
}{h_{m-k}(0)}\mathcal{H}_{m}\left( 0\right)  \\
&\leq &\mathcal{H}_{m}\left( s\right) -\mathcal{H}_{m}\left( 0\right)  \\
&\leq &\sum_{r=1}^{N-1}\sum_{k=0}^{T_{r}}\sum_{l=1}^{Z_{r}(k)}\eta
_{r,N}\left( k,l\right) \frac{\left( h_{m-k}(s)-h_{m-k}(0)\right)
}{h_{m-k}(s)}\mathcal{H}_{m}\left( s\right)  \\
&\leq &W_{N}\left( h_{m}(s)\right) ^{W_{N}-1}\left(
h_{m-T_{N-1}}(s)-h_{m-T_{N-1}}(0)\right) .
\end{eqnarray*}
Hence, on account of the condition $T_{N-1}\leq n^{2/3}$
we conclude that
\begin{eqnarray*}
&&\left( h_{m}(s)-h_{m}(0)\right) \mathbf{E}\left[ W_{N}\left(
h_{m-n^{2/3}}(0)\right) ^{W_{N}}I_{N-1}(n^{2/3})\right]  \\
&\leq &\mathbf{E}\left[ \left( \mathcal{H}_{m}\left( s\right)
-\mathcal{H}_{m}\left( 0\right) \right) I_{N-1}(n^{2/3})\right]  \\
&\leq &\left( h_{m-n^{2/3}}(s)-h_{m-n^{2/3}}(0)\right) \mathbf{E}\left[
W_{N}\left( h_{m}(s)\right) ^{W_{N}-1}I_{N-1}(n^{2/3})\right] .
\end{eqnarray*}

Observe that if the parameter $s$ has form (\ref{defS}),
then for $1\ll k\ll m\sim n\rightarrow \infty $
\begin{equation*}
h_{m}(s)-h_{m}(0)\sim h_{m-n^{2/3}}(s)-h_{m-n^{2/3}}(0)\sim
\frac{k}{b_{N}\lambda n^{2}}
\end{equation*}
in view of Lemma \ref{L_Diff1}, and
 \begin{equation*}
1-h_{m}(s)\sim 1-h_{m}(0)\sim \frac{1}{b_{N}n}
\end{equation*}
by the criticality condition. This, combined with Corollary \ref{C_Cansor}, yields
\begin{eqnarray*}
\Delta (m,k;\lambda ) &\sim &\left( h_{m}(s)-h_{m}(0)\right) \mathbf{E}\left[
W_{N}\left( h_{m}(s)\right) ^{W_{N}}I_{N-1}(n^{2/3})\right]  \\
&\sim &\frac{k}{b_{N}\lambda n^{2}}b_{N}g_{1,N}n^{1-\gamma
_{1}}=\frac{kg_{1,N}}{\lambda n^{1+\gamma _{1}}}.
\end{eqnarray*}
Thus, for $1\ll k\ll m\sim n\rightarrow \infty $
\begin{equation*}
\lim_{n\to\infty}\frac{n^{1+\gamma _{1}}}{k}\Delta (m,k;\lambda )=\frac{g_{1,N}}{\lambda
}.
\end{equation*}
Let now $\lambda$ be a complex variable and
\begin{equation*}
s=\exp \left\{ -\frac{\lambda }{b_{N}k}\right\},\, {Re}\lambda >0.
\end{equation*}
It is not difficult to check that
for ${Re}\lambda \geq \lambda _{0}>0$ and $k\ll m\sim
n\rightarrow \infty $ there exists a constant $C_1(\lambda_0)>0$ such that
 \begin{equation*}
\frac{n^{2}}{k}\left\vert h_{m-n^{2/3}}(s)-h_{m-n^{2/3}}(0)
\right\vert \leq \frac{n^{2}}{k}\left( h_{m-n^{2/3}}(\left\vert s\right\vert
)-h_{m-n^{2/3}}(0)\right) \leq \frac{C_1(\lambda_0)}{\lambda _{0}},\,
\end{equation*}
and, in addition, by Corollary \ref{C_Cansor} there exists a constant $C_2(\lambda_0)>0$  such that
\begin{equation*}
n^{\gamma _{1}-1}\left\vert \mathbf{E}\left[ W_{N}\left( h_{m}(s)\right)
^{W_{N}}I_{N-1}(n^{2/3})\right] \right\vert \leq n^{\gamma_{1}-1}\mathbf{E}\left[ W_{N}\left( h_{m}(\left\vert s\right\vert )\right)
^{W_{N}}I_{N-1}(n^{2/3})\right] \leq C_2(\lambda _{0}).
\end{equation*}
Basing on these estimates and using the inequalities
\begin{eqnarray*}
 \Big|\prod_{j=1}^{J}a_{j}^{r_{j}}-\prod_{j=1}^{J}b_{j}^{r_{j}}\Big|= \Big|\sum_{l=1}^{J}(a_l^{r_l}-b_l^{r_l})\prod_{j=l+1}^{J}a_{j}^{r_{j}}\prod_{k=1}^{l-1}b_{k}^{r_{k}}\Big|\\
\leq \sum_{l=1}^{J}r_{l}|a_{l}-b_{l}||a_{l}|^{r_l-1}\Big|\prod_{j=l+1}^{J}a_{j}^{r_{j}}\prod_{k=1}^{l-1}b_{k}^{r_{k}}\Big|\\
\leq \sum_{l=1}^{J}r_{l}|a_{l}-b_{l}|a^{\sum_{j=1}^Jr_j-1},
\end{eqnarray*}
being valid for $0< |b_{j}|\leq |a_{j}|\leq a\leq 1$
and nonnegative integers $r_{j},$ we conclude that for ${Re}\lambda \geq \lambda _{0}>0$
\begin{eqnarray*}
\frac{n^{1+\gamma _{1}}}{k}\left\vert \Delta (m,k;\lambda )\right\vert
&\leq &\frac{n^{2}}{k}\left(
h_{m-n^{2/3}}(|s|)-h_{m-n^{2/3}}(0)\right) n^{\gamma
_{1}-1}\mathbf{E}\left[ W_{N}\left( h_{m}(|s|)\right)
^{W_{N}-1}I_{N-1}(n^{2/3})\right]\\
&\leq &C_{3}\left( \lambda _{0}\right)
\end{eqnarray*}
for some constant $C_3(\lambda_0)>0.$  This fact allows us, the same as in the proof of Lemma \ref{L_Diff2}, to apply the Montel, Vitali and Weierstrass theorems and to conclude that for $1\ll k\ll n$
\begin{multline*}
\lim_{n\to\infty}\frac{n^{1+\gamma _{1}}}{k^{2}}\mathbf{E}\left[
Z_{N}(m)\exp \left\{ -\lambda \frac{Z_{N}(m)}{b_{N}k}\right\}
I_{N-1}(n^{2/3})\right]  \\
=-\frac{\partial }{\partial \lambda }\lim_{n\to\infty}\frac{b_{N}n^{1+\gamma _{1}}}{k}\Delta (m,k;\lambda )=\frac{b_{N}g_{1,N}}{\lambda
^{2}}.
\end{multline*}

The lemma is proved.

\section{Proofs of Theorems \protect\ref{T_death} and  \protect\ref{T_deathFin}}

\ \textbf{Proof of Theorem \ref{T_death}}. By virtue of
(\ref{NoPrevious}) it is sufficient to show that, for $n\sim
m\gg k\rightarrow \infty $ and $k=n-m$
\begin{equation*}
\mathbf{E}\left[ \exp \left\{ -\lambda _{N}\frac{Z_{N}(m)}{b_{N}k}\right\}
I_{N-1}(n^{2/3})\Big|\,T_{N}=n\right] \rightarrow \frac{1}{\left( \lambda
_{N}+1\right) ^{2}}.
\end{equation*}
Put
 \begin{equation*}
s=\exp \left\{ -\frac{\lambda _{N}}{b_{N}k}\right\} .
\end{equation*}
Clearly, for $m>n^{2/3}$
\begin{equation*}
\mathbf{E}\left[ s^{Z_{N}(m)}I_{N-1}(n^{2/3});\,T_{N}=n\right]
=\mathbf{E}\left[ \left( \left( sh_{k+1}(0)\right) ^{Z_{N}(m)}-\left( sh_{k}(0)\right)
^{Z_{N}(m)}\right) I_{N-1}(n^{2/3})\right]
\end{equation*}
and, therefore, as $k\to\infty$
\begin{multline*}
s\left( h_{k+1}(0)-h_{k}(0)\right) \mathbf{E}\left[ Z_{N}(m)\left(
sh_{k}(0)\right) ^{Z_{N}(m)}I_{N-1}(n^{2/3})\right]  \\
\leq \mathbf{E}\left[ s^{Z_{N}(m)}I_{N-1}(n^{2/3});\,T_{N}=n\right]  \\
\leq s\left( h_{k+1}(0)-h_{k}(0)\right) \mathbf{E}\left[ Z_{N}(m)\left(
sh_{k+1}(0)\right) ^{Z_{N}(m)-1}I_{N-1}(n^{2/3})\right] .
\end{multline*}
By the local limit theorem for the critical Galton-Watson processes
\begin{equation*}
h_{k+1}(0)-h_{k}(0)\sim \frac{1}{b_{N}k^{2}}, \quad k\to\infty.
\end{equation*}
Further,
\begin{equation*}
sh_{k+1}(0)=\exp \left\{ -\frac{\lambda _{N}}{b_{N}k}+\log
h_{k+1}(0)\right\} =\exp \left\{ -\frac{\lambda _{N}+1}{b_{N}k}\left(
1+o(1)\right) \right\} .
\end{equation*}
Using Lemma \ref{L_Diff3} we conclude that for $n\sim
m\gg k\rightarrow \infty $ and $k=n-m$
\begin{equation*}
\mathbf{E}\left[ Z_{N}(m)\exp \left\{ -\left( \lambda _{N}+1\right) \left(
1+o(1)\right) \frac{Z_{N}(m)}{b_{N}k}\right\} I_{N-1}(n^{2/3})\right] \sim
\frac{k^{2}b_Ng_{1,N}}{n^{1+\gamma _{1}}\left( \lambda _{N}+1\right) ^{2}}.
\end{equation*}
Since, as $n\to\infty$
\begin{equation*}
\mathbf{P}\left( T_{N}=n\right) \sim \frac{g_{1,N}}{n^{1+\gamma _{1}}}
\end{equation*}
in view of (\ref{LocAsymp}), we have
\begin{eqnarray*}
\mathbf{E}\left[ s^{Z_{N}(m)}I_{N-1}(n^{2/3});\,T_{N}=n\right]  &\sim
&\frac{1}{b_{N}k^{2}}\frac{k^{2}b_{N}g_{1,N}}{n^{1+\gamma _{1}}\left( \lambda
_{N}+1\right) ^{2}} \\
&=&\frac{g_{1,N}}{n^{1+\gamma _{1}}\left( \lambda _{N}+1\right) ^{2}}\sim
\mathbf{P}\left( T_{N}=n\right) \frac{1}{\left( \lambda _{N}+1\right) ^{2}},
\end{eqnarray*}
as required.

\textbf{Proof of Theorem \ref{T_deathFin}.} The same as in the proof of Theorem \ref{T_death}, it is sufficient to show that, for a fixed
 $k=n-m$ and $n\rightarrow \infty $
\begin{equation*}
\mathbf{E}\left[ s^{Z_{N}(m)}I_{N-1}(n^{2/3})\Big|\,T_{N}=n\right]
\rightarrow s\big(U\left( sh_{k+1}(0)\right) -U\left( sh_{k}(0)\right)\big).
\end{equation*}
Put
\begin{equation*}
k^{\ast}=\left[ \sqrt{n}\right],\quad    m^{\ast }=m-\left[ \sqrt{n}\right] =m-k^{\ast}
\end{equation*}
and write the representation
\begin{eqnarray*}
&&\mathbf{E}\left[ s^{Z_{N}(m)}I_{N-1}(n^{2/3});\,T_{N}=n\right]\\
&=&\mathbf{E}\left[ \left( \left( h_{k^{\ast}}\left( sh_{k+1}(0)\right)
\right) ^{Z_{N}(m^{\ast })}-\left( h_{k^{\ast}}\left( sh_{k}(0)\right)
\right) ^{Z_{N}(m^{\ast })}\right) I_{N-1}(n^{2/3})\right] .
\end{eqnarray*}
Similarly to the estimates used in the proof of Theorem \ref{T_death}, we have
\begin{eqnarray*}
&&s\left( h_{k^{\ast}}\left( sh_{k+1}(0)\right) -h_{k^{\ast}}\left(
sh_{k}(0)\right) \right) \mathbf{E}\left[ Z_{N}(m^{\ast })\left(
h_{k^{\ast}}\left( sh_{k}(0)\right) \right) ^{Z_{N}(m^{\ast
})}I_{N-1}(n^{2/3})\right]  \\
&\leq &\mathbf{E}\left[ s^{Z_{N}(m)}I_{N-1}(n^{2/3});\,T_{N}=n\right]  \\
&\leq &s\left( h_{k^{\ast}}\left( sh_{k+1}(0)\right) -h_{k^{\ast}}\left( sh_{k}(0)\right) \right) \mathbf{E}\left[ Z_{N}(m^{\ast })\left(
h_{k^{\ast}}\left( sh_{k+1}(0)\right) \right) ^{Z_{N}(m^{\ast
})-1}I_{N-1}(n^{2/3})\right] .
\end{eqnarray*}
According to (\ref{DefHarm})
\begin{equation*}
\lim_{n\rightarrow \infty }\left( k^{\ast}\right) ^{2}\left(
h_{k^{\ast}}\left( sh_{k+1}(0)\right) -h_{k^{\ast}}\left(
sh_{k}(0)\right) \right)=b_{N}^{-1}\left( U(sh_{k+1}(0)\right) -U(sh_{k}(0))).
\end{equation*}
Since
\begin{equation*}
1-h_{k^{\ast}}\left( sh_{k+1}(0)\right) \sim 1-h_{k^{\ast}}\left(
0\right) \sim \frac{1}{b_{N}k^{\ast}}
\end{equation*}
as $n\rightarrow \infty, $ and
$k^{\ast}\ll m^{\ast }\sim n,$ then, according to Lemma \ref{L_Diff3}
\begin{equation*}
\lim_{n\to\infty}\frac{n^{1+\gamma _{1}}}{\left(k^{\ast}\right)^{2}}\mathbf{E}\left[ Z_{N}(m^{\ast })\exp \left\{
-\frac{Z_{N}(m^{\ast })}{b_{N}k^{\ast}}\right\} I_{N-1}(n^{2/3})\right].
=b_Ng_{1,N}
\end{equation*}
Thus, for any fixed $k=n-m$ and $n\rightarrow \infty$
\begin{eqnarray*}
\mathbf{E}\left[ s^{Z_{N}(m)}I_{N-1}(n^{2/3});\,T_{N}=n\right]  &\sim
&\frac{s\left( U(sh_{k+1}(0)\right) -U(sh_{k}(0)))}{b_{N}(k^{\ast})^{2}}\frac{b_Ng_{1,N}(k^{\ast})^{2}}{n^{1+\gamma _{1}}} \\
&\sim &s\left( U(sh_{k+1}(0)\right) -U(sh_{k}(0)))\mathbf{P}\left(
T_{N}=n\right),
\end{eqnarray*}
as required.

Theorem \ref{T_deathFin} is proved.

\end{document}